\documentclass[12pt]{iopart}

\usepackage{iopams}
\usepackage{xspace}
\usepackage{graphicx}

\begin{document}

\title{Constructing a polynomial whose nodal set is the three-twist knot $5_2$}

\author{Mark~R~Dennis and Benjamin~Bode}

\address{H H Wills Physics Laboratory, University of Bristol, Bristol BS8 1TL, UK}

\begin{abstract}
We describe a procedure that creates an explicit complex-valued polynomial function of three-dimensional space, whose nodal lines are the three-twist knot $5_2$.
The construction generalizes a similar approach for lemniscate knots: a braid representation is engineered from finite Fourier series and then considered as the nodal set of a certain complex polynomial which depends on an additional parameter.
For sufficiently small values of this parameter, the nodal lines form the three-twist knot.
Further mathematical properties of this map are explored, including the relationship of the phase critical points with the Morse-Novikov number, which is nonzero as this knot is not fibred.
We also find analogous functions for other knots with six crossings. 
The particular function we find, and the general procedure, should be useful for designing knotted fields of particular knot types in various physical systems.  
\end{abstract}


\section{Introduction}\label{sec:int}

It is often of interest to find an explicit function $f_{K}:S^3 \longrightarrow \mathbb{C}$ whose preimage of zero $f_K^{-1}(0)$ has the topology of some specified knot or link $K$. 
Such functions are known to exist; however, the form of such an $f_K$ is only known explicitly for a number of special cases of $K$.
Knotted functions are not only of mathematical interest, but are important for physical applications; in recent years such functions $f_K$ have been the basis of explicit knotted structures embedded in physical fields, even admitting experimental realization.

In 1928, Brauner found a simple set of functions with zeros in the form of torus knots and links \cite{brauner} (in an approach popularised by Milnor \cite{milnor}), and later Perron \cite{perron}, Rudolph \cite{rudolph} and Dennis et al.~\cite{isolated} found examples of $f_K$ with $K$ as the figure-eight knot and its higher-periodic generalizations: the link of three borromean rings $L6a4$ (also referred to as $6_2^3$), the knots $8_{18}$, $10_{123}$, etc. 
The standard knot table \cite{rolfsen,knotatlas,knotinfo} begins with the trefoil knot $3_1$ (a torus knot), followed by the figure-eight knot $4_1$ and the cinquefoil knot $5_1$ (another torus knot) and then the three-twist knot $5_2$.
Therefore the simplest knot $K$ for which there is no known $f_{K}$ is $5_2$. 
Here, we construct such an $f_{5_2}$ explicitly (Equations (\ref{eq:fa}), (\ref{eq:f5_2}) below, whose nodal knot is plotted in Figure \ref{fig:5_2}), exemplifying a general method which explicitly generates $f_K:S^3 \longrightarrow \mathbb{C}$ with a nodal set isotopic to any knot or link $K$; we describe the general proofs of this method in \cite{proof}.

The knot $5_2$ is distinguished as the simplest non-fibred knot; for a knot $K$ which is not fibred, the complement $S^3 \setminus K$ cannot be represented as a fibre bundle over the base space $S^1$ \cite{rolfsen}.
When $K$ is fibred, a function $f_K$ might provide an explicit fibration provided the derivative $\nabla(\arg f_K) \neq 0$ everywhere, where the fibres are the surfaces of constant argument $\arg f_K$.
In this case, the surfaces of constant argument (phase) are Seifert surfaces, whose boundary is the knot.
Since $5_2$ is not fibred, there must be critical points of $\arg (f_{5_2}(\bullet))$ in $S^3$, in the sense that $\nabla \arg f_{5_2} = 0$ at them.
The level set of constant $\arg f_{5_2}$ at these points is not a 2-surface -- our explicit $f_{5_2}$ therefore allows us to explore an example in detail of how knots can fail to be fibred.

Brauner's knot functions have a very simple form, in terms of complex coordinates $u$, $v$, of the unit 3-sphere $\{ (u,v) \in \mathbb{C}^2: |u|^2 + |v|^2 = 1\}$. 
If $T_{p,q}$ denotes the $p,q$ torus knot or link with $p,q \in \mathbb{N}$, the complex polynomial $f_{T_{p,q}}(u,v) = u^p - v^q$ has a zero set $f_{T_{p,q}}^{-1}(0) \simeq T_{p,q} \subseteq S^3$, as explained in detail at the beginning of \cite{milnor}.
Explicitly, the Hopf link has $(p,q) = (2,2)$ and the trefoil knot has $(p,q) = (2,3)$.
These functions are also holomorphic in $u$ and $v$, and it is straightforward to generalize them to include cables of torus knots \cite{brauner}.
Rudolph's map for the figure-eight knot $f_{4_1} = u^3 - 3 v^2 \overline{v}^2 (1 + v^2 - \overline{v}^2) u - 2(v^2 + \overline{v}^2)$, on the other hand, is of the kind $f_{4_1}:\mathbb{C}\times \mathbb{R}^2 \longrightarrow \mathbb{C}$, which we call \emph{semiholomorphic}: $f$ is holomorphic in $u$, but depends on both $v$ and its conjugate $\overline{v}$.

The functions introduced by Brauner, Perron and Rudolph were not motivated by finding explicit knotted functions on a 3-dimensional manifold like $S^3$, but rather to understand the properties of isolated singularities of more general maps $f:\mathbb{C}^2 \longrightarrow \mathbb{C}$ \cite{milnor}.
Both Brauner's maps $f_{T_{p,q}}(u,v)$ and Rudolph's $f_{4_1}(u,v,\overline{v})$ have a singularity (i.e.~a critical point where $\nabla f = 0$) at the origin in $\mathbb{C}^2$ and $\mathbb{R}^4$ respectively, such that the intersection of the zero sets $f^{-1}_K(0)$ with any 3-sphere of sufficiently small radius gives $K$. 
Akbulut and King \cite{ak} proved non-constructively that a polynomial map exists of the form $\mathbb{R}^4 \longrightarrow \mathbb{C}$ which has a weakly isolated singularity for any knot or link.

Our recent generalisation \cite{lemniscate} of the approach of \cite{isolated} for a different choice of $f_{4_1}$ to the family of lemniscate knots (described below) was similarly semiholomorphic; the knots here only occur on restricting the polynomial in $u,v,\overline{v}$ to the unit 3-sphere (and other radii close to $1$) and not necessarily to 3-spheres of other radii, and at best have weakly isolated singularities in the sense of \cite{ak}.
The lemniscate knot construction and its generalizations \cite{lemniscate} are, however, highly suited as candidates for embedding knots into physical fields in a range of physical systems, generalizing Brauner's maps for torus knots. 
Of course, physical systems are more usually defined in three-dimensional cartesian space than $S^3$, but the function $f_K(u,v)$ can be generalised to $f_K(x,y,z)$ using standard stereographic coordinates
\begin{equation} 
   \left. \begin{array}{cll} 
   u & = \displaystyle \frac{x^2+y^2+z^2 - 1 + 2 \rmi\, z}{x^2+y^2+z^2 + 1} & = \displaystyle \frac{R^2+z^2 - 1 + 2 \rmi\, z}{R^2+z^2 + 1}  \\
   v & = \displaystyle \frac{2(x+\rmi\, y)}{x^2+y^2+z^2 + 1} & \displaystyle = \frac{2 R \exp(\rmi\, \phi)}{R^2+z^2 + 1}, \end{array} \right\}
   \label{eq:stereo}
\end{equation}
where $x,y,z$ are cartesian coordinates, and $R,\phi, z$ are cylindrical polar coordinates which we will often find more useful in discussing the construction.

Any explicit representation of a knot must start with a way of encoding the knot's conformation uniquely; our construction of $f_{5_2}$, as with the other knotted field functions discussed is based on a braid representation, of which the knot is the closure. 
The torus knot $T_{p,q}$ is the closure of a $p$-strand helix which rotates by $q/p$ turns before closing; a figure-eight knot is the closure of a 3-strand pigtail braid after it undergoes two periods (and hence four crossings) of its characteristic over-under pattern.
In terms of the Artin braid group \cite{kt}, a braid of $s$ strands has $s-1$ generators $\sigma_k$, $k,1,\dots,s-1$, where $\sigma_k$ represents an overcrossing between the $k$th and $k+1$th strands, and $\sigma_k \sigma_{k+1} \sigma_k = \sigma_{k+1} \sigma_k \sigma_{k+1}$ for each $k$, and other generators commute. 
The torus knot $T_{s,q}$ is therefore represented by the braid word $(\sigma_1\sigma_2\cdots\sigma_{s-1})^q$; the figure-eight knot by $(\sigma_1\sigma_2^{-1})^2$ (and higher powers $(\sigma_1\sigma_2^{-1})^n$ for higher periods of the knot, such as $n = 3$ for $L6a4$, $n = 4$ for $8_{18}$ and $n = 5$ for $10_{123}$).

The lemniscate knots are a large family of knots and links including the torus knots and figure-eight knot, given by closures of braids that can be parametrized trigonometrically as follows.
Positive integers $s, \ell, r$ (with $s$, $\ell$ coprime) are chosen so $\ell$ determines a $(1,\ell)$ Lissajous figure $T_L$ in the form of a generalized lemniscate with respect to $t \in [0,2\pi s)$,
\begin{equation}
   T_L: (X_L(t), Y_L(t)) = \left( \cos\left(\frac{t}{s}\right), \frac{1}{\ell} \sin\left(\frac{\ell t}{s}\right)\right).
   \label{eq:genlem}
\end{equation} 
The factor $1/\ell$ is chosen so when $\ell = 2$ this is a lemniscate of Gerono. 
The $(s,\ell,r)$-\emph{lemniscate knot} $L = L(s,\ell,r)$ is the closure of the braid 
\begin{equation}
  \mathcal{B}_L = \bigcup_{j=0}^{s-1}\left(X_L(rt + 2\pi j), Y_L(rt + 2\pi j),t\right), \qquad t \in [0,2\pi).
  \label{eq:lembraid}
\end{equation}
As the braid length parameter $t$ increases, the $s$ strands of the braid execute the path $T_L$ in the transverse plane of the braid, perpendicular to the direction of increase of $t$.
The strands are distributed equally around the curve $T_L$ (with respect to its parametrization by $t$), and $s$ and $\ell$ must be coprime to ensure the strands never intersect.
As $t$ increases from $0$ to $2\pi$, each strand makes its way $r/s$ around $T_L$; as $t \to t+2\pi/r$, the strands are cyclically permuted, $1 \to 2 \to 3 \to \cdots \to s \to s+1 = 1$. 
One can read off the braid word by drawing the two-dimensional braid diagram in the $(t,X_L)$ plane, where crossing signs are found from $Y_L$; by (\ref{eq:genlem}), when $t = 2\pi k/r$, $k = 1, \ldots r-1$ all crossings labelled by an odd generator occur, with all even generators at at $t = \pi (1+2k)/r$.
The word therefore has the form $(\sigma_1^{\epsilon_1} \sigma_3^{\epsilon_3} \cdots \sigma_2^{\epsilon_2} \sigma_4^{\epsilon_4} \cdots)^r$, where the signs $\varepsilon_j = \pm 1$ for $j = 1, \ldots, s$ are determined by the sign of $Y_L(t)$ at $t$ evaluated at integer multiples of $\pi$, corresponding to the crossings.
It is therefore clear that the torus knot/link $T_{p,q}$ is represented in this way by $(s,\ell,r) = (p,1,q)$, and $4_1$ by $(s,\ell,r) = (3,2,2)$ (higher-period generalizations have higher values of $r$).

The rest of the construction follows by encoding the braid into a family of polynomials $p_t(u)$, which have, for each $t$, roots at $u = a X_L(rt + 2\pi j) + \rmi b Y_L(rt + 2\pi j)$ for $j = 1, \ldots s$ and $a$ and $b$ real nonzero scaling factors.
Multiplying the polynomial out and using simple arithmetic of roots of unity, it can be shown that $p_t(u)$ is a polynomial in $u, \rme^{\rmi t}, \rme^{-\rmi t}$ (i.e.~no fractional powers of the exponential). 
One finds a function $f(u,v,\overline{v})$ by replacing $\rme^{\rmi t} \to v$ and $\rme^{-\rmi t} \to \overline{v}$ in this polynomial expression; for sufficiently small $a$ and $b$, this function is indeed $f_L$, and is automatically a semiholomorphic polynomial.
The trigonometry of the braid function has thus been replaced by algebra in the polynomial.

Of course, not all knots are in the lemniscate family (in fact all lemniscate knots are fibred since the lemniscate braid representation is homogeneous \cite{stallings}): more complicated functions are needed to realize braid representations of other knots.
A convenient---although by no means unique---braid representative for a knot or link is its minimal braid representation \cite{gittings}, defined as the shortest braid word closing to the knot (which may be longer than the knot's crossing number, and also is not necessarily unique).
For $5_2$, this is given by the word of length $6$ with two generators, $w \equiv \sigma_1^{-1} \sigma_2\sigma_1^3 \sigma_2$; the braid word $w$ is not homogeneous -- both $\sigma_1$ and $\sigma_1^{-1}$ occur \cite{stallings}.

In order to find $f_{5_2}$, it is natural to generalize the lemniscate knot construction; we must find a pair of trigonometric functions $X(t), Y(t)$ which encode the minimal braid word $w$ in the same way that $X_L(t),Y_L(t)$ encode the lemniscate knot's braid word.
After encoding the braid into a polynomial in $u,v,\overline{v}$, we must find a suitable choice for $a$ and $b$.
We will go through these steps explicitly in detail, to illustrate various features of the construction. 
The braid functions are found by inspection, using simple finite Fourier series instead of (\ref{eq:genlem}).
The behaviour of the function on different scaling factors is described in some numerical detail, which depends on critical points of the absolute value $|f_{5_2}|$. 
These are distinct from the critical points of phase (i.e.~argument) $\arg f_{5_2}$, whose existence is related to the failure of $5_2$ to be fibred, and these also can easily be located numerically.
We will also describe how the $f_{5_2}$ we construct cannot easily be generalized to a function with an isolated singularity.
Finally, we consider similar construction of functions with nodal knots for other knots of six and fewer crossings. 
Before proceeding, however, we will briefly survey the various physical systems in which knotted fields such as $f_{5_2}$ might be of relevance.

There is much interest, going back to Kelvin \cite{kelvin}, in studying the evolution of vortex knots in fluids (either classical or superfluid) \cite{kleckner,hall}, and this construction provides natural phase functions whose gradient describes the velocity pattern around a knotted vortex.
The construction in \cite{kedia} gives a precise prescription to create a classical knotted flow field based on a complex scalar function $f_K$ with a nodal knot; this has the additional feature that, since the underlying complex function is semiholomorphic, the helicity of the flow can be chosen to be any integer multiple of the number of strands in the construction of the knot.
Electromagnetic fields with knotted magnetic field lines can also be constructed by an analogous construction \cite{kedia1}.

Vortex lines, as nodal lines of a complex scalar field, also occur in quantum eigenfunctions, and can be designed to be knotted and linked \cite{berry}; so far Brauner's functions, or others, have not been embedded into stationary eigenfunctions of a quantum system.
It is known that knotted nodal lines exist in the nodal sets of generic eigenfunctions in simple quantum systems such as the three-dimensional harmonic oscillator or hydrogen atom \cite{daniel1,daniel2,taylor}, but these knots are components of densely linked tangles, and not the isolated nodal sets considered here.
Propagating, coherent optical fields satisfy the paraxial wave equation, which is the $2+1$ Schr\"odinger's equation (the third propagation direction corresponds to time); with the $z = 0$ section of a Brauner function or a figure-eight generalization, the propagating field also has a nodal knot as realised experimentally in \cite{isolated}.
Complex scalar functions are also useful in describing reaction-diffusion systems (such as in chemical waves) in which the zeros are organizing centres, and knots have been studied in such systems mathematically \cite{winfreestrogatz,winfreesutcliffe,sm}.

Physical systems described by other target spaces can have knotted fields based on nodal sets of complex scalar functions.
One such example is topological solitons (described by maps from the three-sphere to the two-sphere), in which complex function-based rational maps have been successful as ans\"atze for low-lying energy states in the Skyrme-Faddeev model \cite{sutcliffe}.
Stationary configurations of energy are known for this model to take the form of torus knots, cable knots and links \cite{sutcliffe,jennings}, but at present no other knots have been found; an initial ansatz with the topology of some other knot (such as $5_2$) can be made following the description in \cite{lemniscate}.

Liquid crystals are another class of systems area in which knotted defects can be created and controlled \cite{ma,copar}.
A particularly interesting possibility for knotted defects is in smectic liquid crystals, whose low-energy configurations are in layers of material naturally represented by the complex argument of a complex scalar field, and in which realisations of fibred knots \cite{km} seem most natural. 
Understanding how smooth phase fields break down for non-fibred knots is then of particular interest in these systems.

\section{Finding a suitable trigonometric braid}\label{sec:braid}

\begin{figure}[htbp]
\begin{center}
\includegraphics[width=14cm]{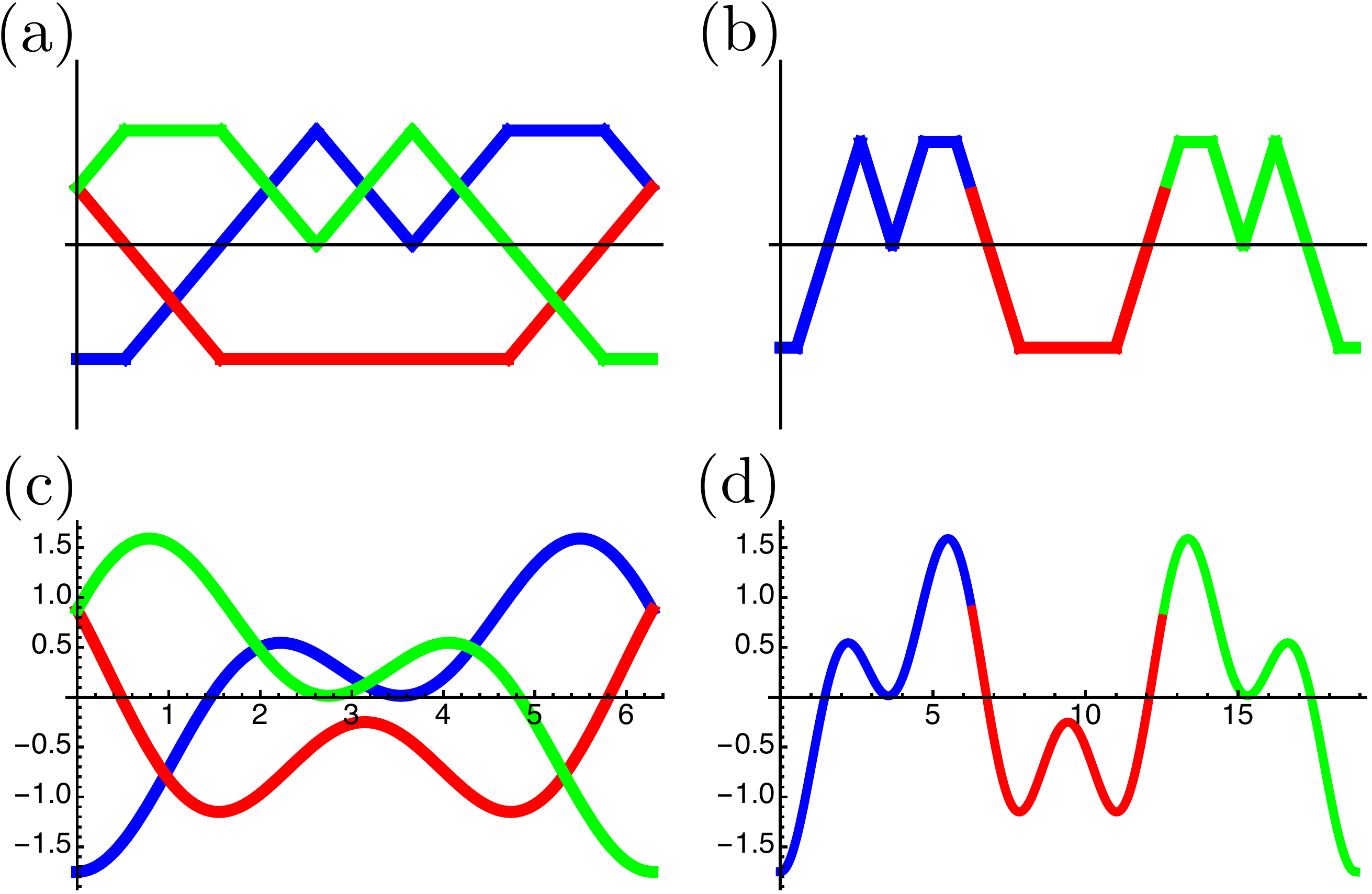} 
\caption{
   Constructing $X(t)$ from the unsigned braid diagram of $5_2$.
   (a) Piecewise linear representation of the unsigned word $w^{\bullet}$, with three strands.
   Crossing $\sigma_1$ occurs above crossing $\sigma_2$, and strands are labelled increasing downwards at $t=0$: $1$ (green), $2$ (red), $3$ (blue).
   (b) Unfolded braid diagram made by placing the strands in (a) end-to-end.
   (c) Plot of $X(t),X(t+2\pi),X(t+4\pi)$ from (\ref{eq:Xchoice}), giving the same unsigned braid diagram as (a) for $t \in [0,2\pi]$.
   (d) Unfolded $X(t)$ from (\ref{eq:Xchoice}) plotted from $t \in [0,6\pi]$.
}
\label{fig:interpolation}
\end{center}
\end{figure}

We seek explicit functions $X(t), Y(t)$ for the three-twist knot $5_2$ such that the braid
\begin{equation}
   \mathcal{B}_{5_2} = \bigcup_{j=0}^{2}\left(X(t + 2\pi j), Y(t + 2\pi j),t\right), \qquad t \in [0,2\pi)
   \label{eq:B52}
\end{equation}
is a representation of its minimal braid word $w = \sigma_1^{-1} \sigma_2\sigma_1^3 \sigma_2$.
As discussed above, if the two-dimensional braid diagram is represented in the $(t,X)$ plane, the function $X(t)$ encodes only the sequence of generators and not their signs, $w^{\bullet} \equiv {\sigma_1^{\bullet}} \sigma_2^{\bullet} {\sigma_1^{\bullet}}^3 \sigma_2^{\bullet}$, where $\bullet$ on a generator denotes an intersection of strands whose crossing sign $\pm 1$ is undetermined.
It is convenient to place the first crossing $t = t_1$ at $t_1 = 0$.
This 3-strand braid (without crossing signs) is represented as a piecewise-linear graph with respect to $t \in [0,2\pi]$ in Figure \ref{fig:interpolation} (a), where strands and crossings are labelled increasing from the top left downwards. 
Upon closure, the three strands are joined end-to-end to form a single knotted component; the single composite strand in $(t,X)$ is a function with $t \in [0,6\pi]$ as shown in Figure \ref{fig:interpolation} (b).
This periodic function is in fact symmetric about $t = 0$ and $t = 3\pi$. 
The desired function $X(t)$ should be a trigonometric function---conveniently a finite Fourier series---which approximates this piecewise linear form such that, $X(t), X(t+2\pi)$ and $X(t+4\pi)$ give the right pattern of intersections for $w^{\bullet}$.

Over the $6\pi$ period of the function in Figure \ref{fig:interpolation} (b), there is a larger amplitude modulation of period $2$, say $-\cos(2t/3)$.
This is accompanied by a smaller amplitude modulation of higher spatial frequency not rationally related to the larger wave.
The two sinusoidal components should have most interference around $t = \pi$ and $5\pi$, so spatial frequency $5/3$ appears suitable. 
Therefore a two-term Fourier series, with appropriately chosen amplitude, is sufficient to define $X(t)$, 
\begin{equation}
   \textstyle X(t) = - \cos(\frac{2}{3}t) - \frac{3}{4} \cos(\frac{5}{3}t), 
   \label{eq:Xchoice}
\end{equation}
which is plotted in Figure \ref{fig:interpolation} (c) and (d) analogous to the piecewise linear representation of (a) and (b).
The three curves $X(t), X(t+2\pi)$ and $X(t+4\pi)$ intersect six times for $t \in [0,2\pi)$, $t_1 = 0$, $t_6 = 2\pi-t_2$ (between strands $1$ (green) and $2$ (red)), $t_3 = 1.98\ldots$, $t_4 = \pi$, $t_5 = 2\pi-t_3$ (between strands $2$ and $3$ (blue)), $t_2 = 0.96\ldots$ (between strands $3$ and $1$), where decimal forms have been given for crossings whose explicit expression can be given in terms of arctangents of the roots of the polynomial $9 - 244 z^2 + 798 z^4 - 1092 z^6 + 161 z^8$. 
The process of finding a suitable $X(t)$ requires following the strands, rather than braid algebra which depends simply on the crossings.

\begin{figure}[htbp]
\begin{center}
\includegraphics[width=12cm]{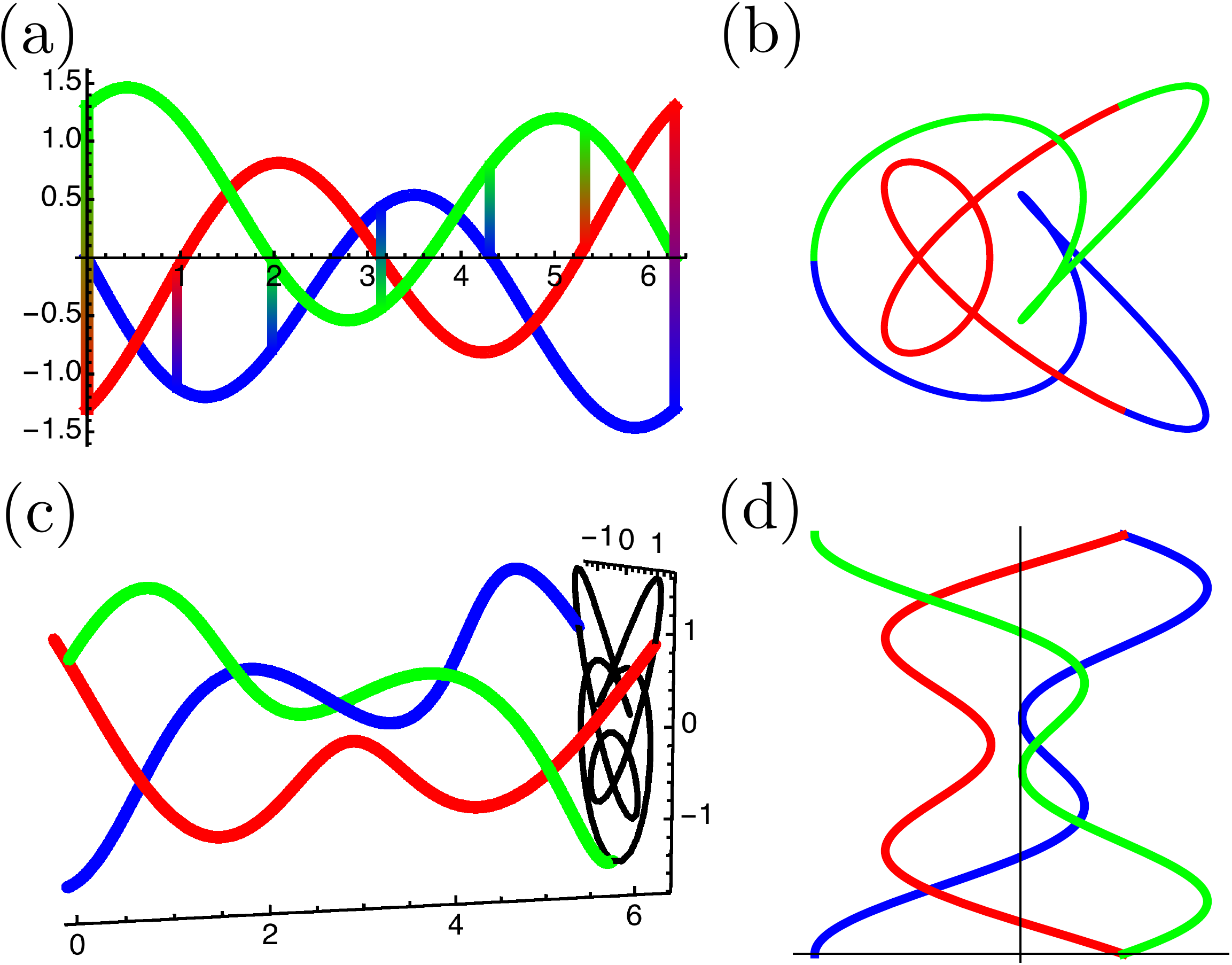}
\caption{Trigonometric braid for $5_2$.
(a) Plots of $Y(t), Y(t+2\pi), Y(t+4\pi)$ showing positions of crossings $t_i$ from $X(t)$.  
(b) The parametric trajectory $T_{5_2}$ given by $(X(t),Y(t))$ for $0 \le t < 6\pi$, colour coded to show $t\in[0,2\pi)$ (strand 3, blue), $t \in [2\pi, 4\pi)$ (strand $2$, red), $t \in [4\pi,6\pi)$ (strand $1$, green).
(c) The braid $\mathcal{B}_{5_2}$, with the trajectory as its transverse projection. 
(d) Trigonometric braid diagram from $X(t)$ projecting to trajectory in (b).
}
\label{fig:braid}
\end{center}
\end{figure}

The function $Y(t)$ is chosen in a similar way; it should be a simple Fourier series correctly encoding the crossing signs.
The sequence of crossing signs in $w$ is $-+++++$, and in terms of the strands, this means: at $t_1$, strand $1$ (green) passes over $2$ (red); at $t_2$, strand $2$ over $3$ (blue); $t_3$, $1$ over $3$; $t_4$, $3$ over $1$; $t_5$, $1$ over $3$; $t_6$, $1$ over $2$. 
In fact, by inspection this is almost satisfied by $-\sin(4t/3)$, with only the crossings at $t_3$ and $t_5$ incorrect. 
As before, adding one other term corrects this, and a suitable choice gives
\begin{equation}
   \textstyle Y(t) = -\sin(\frac{4}{3}t) - \frac{1}{2} \cos(\frac{1}{3}t),
   \label{eq:Ychoice}
\end{equation}
which is antisymmetric about $t = 0$ and $3\pi$, and evidently satisfies the necessary conditions. 
This function is plotted in Figure \ref{fig:braid} (a).

All of the information about the trigonometric braid for $5_2$ is therefore in the \emph{trajectory curve} $T_{5_2}$, given by $(X(t),Y(t))$ for $t \in [0,6\pi)$, shown in Figure \ref{fig:braid} (b) which plays the same role as the $T_L$ Lissajous figure for lemniscate knots.
It has a 2-fold reflection symmetry due to the symmetry-antisymmetry of $X$ and $Y$.
It is the transverse trajectory of the braid $\mathcal{B}_{5_2}$ plotted in Figure \ref{fig:braid} (c), which correctly gives the braid word $w$.
Finding this trigonometric braid has solved the main problem in constructing $f_{5_2}$.

Different projection directions of the braid $\mathcal{B}_{5_2}$ give different braid words that must also close to $5_2$.
In particular, projecting in the $(t,Y)$ plane gives the 8-crossing palindromic braid word $\sigma_1 \sigma_2^{-1} \sigma_1 \sigma_2 \sigma_1 \sigma_2 \sigma_1^{-1} \sigma_2$ (evident from Figure \ref{fig:braid} (a)).

It is not difficult to find alternative functions $X(t)$ and $Y(t)$ which also satisfies the requirements of crossings (e.g.~for the same $X(t)$, $Y(t) = - \sin(\frac{4}{3}t) + \frac{8}{5}\sin(\frac{7}{3}t)$); we believe our choice is close to being as simple as possible, with a minimal number of Fourier terms, and the simplest rational frequencies and coefficients.
The overall magnitudes of the maximum of $X(t)$ and $Y(t)$ will be important later; we observe that $\max_{t\in[0,6\pi)}|(X(t), Y(t))| =  
2.098\ldots$ at $t = t_{\mathrm{max}} \equiv 5.535\ldots$ and $t = 6\pi - t_{\mathrm{max}}$.

\section{From trigonometric braid to nodal knot}\label{sec:embed}

Having constructed an appropriate braid function, we now use it to create functions with the braid topology encoded into their nodal pattern -- this will eventually lead to the desired $f_{5_2}:S^3\longrightarrow \mathbb{C}$ with $f^{-1}_{5_2}(0) \simeq 5_2$.
The reasoning behind the approach (which follows the logic described in \cite{lemniscate} closely, and also \cite{rudolph}) is as follows. 
First a function is constructed from the three-dimensional space in which the braid is embedded in $\mathbb{C}$, which has zeros along the braid $\mathcal{B}_{5_2}$ and is a polynomial, as we will see.
We wish to embed the closed braid zeros in a function whose domain is $S^3$.
We do this by generalizing the braid polynomial to a suitable function in a 4-dimensional space, whose intersection with a restricted 3-dimensional subspace gives the nodal braid, and on restricting to the unit 3-sphere gives the correct nodal knot.
Even with the correct function this is not automatic, and relies in fact on a single tuning parameter, whose values to give $5_2$ in practice are determined numerically.

We construct the braid polynomial as a map $p_t:\mathbb{C}\times S^1 \longrightarrow \mathbb{C}$ whose nodal set is $\mathcal{B}_{5_2}$, so the transverse plane of the braid is the complex Argand plane (in the variable $u$), and $t \in [0,2\pi)$ is the cyclic braid parameter. 
This is clearly 
\begin{eqnarray}
\fl   p_{t,a}(u) & = & \prod_{j = 1}^3 (u - a[X(t + 2\pi j) + \rmi Y(t+2\pi j)]) \nonumber \\
\fl    & = & 256 u^3 - 12 a^2 u \{5 + 8 \cos(t) + 4\rmi\,[4 \sin(t) + 11 \sin(2 t) + 6 \sin(3 t)]\} \nonumber \\
\fl    & & + a^3 \left\{ 192 + 372 \cos(t) + 256 \cos(2 t) + 144 \cos(3 t) + 108 \cos(4 t) + 27 \cos(5 t) \right.\nonumber \\ 
\fl    & & \left. - 2 \rmi\, [188 \sin(t) + 102 \sin(2 t) - 21 \sin(3 t) - 32 \sin(4 t)]\right\}. \label{eq:pta}
\end{eqnarray}
Here, $a > 0$ has been introduced to rescale the braid's geometric width, i.e.~the trajectory in the complex $u$ plane; it is the tuning parameter we will need in the following.
With this map, the braid is now implicitly closed and $t \in S^1$.
For each value of $t$, $p_t(u)$ is simply a complex cubic polynomial, whose three roots sweep out $\mathcal{B}_{5_2}$ as $t$ increases.

The polynomial (\ref{eq:pta}) is now generalized to a map $f_a(u,v):\mathbb{C}\times\mathbb{R}^2 \longrightarrow \mathbb{C}$, such that the restriction of $f_a$ to $\mathbb{C}\times S^1 \subseteq \mathbb{C}\times\mathbb{R}^2$ is $p_{t,a}$, and the restriction of $f_a$ to the three-sphere $S^3 \subseteq \mathbb{C}\times\mathbb{R}^2$ gives our desired function $f_{5_2}$, for an appropriate choice of $a$.
This generalization is made by considering $t$ as the complex argument of a new complex variable $v$, so $t = \arg v$; since $\cos(t)$ and $\sin(t)$ appear in (\ref{eq:pta}) on unequal footing, the generalized function $f_a(u,v,\overline{v})$ depends both on $v$ and its conjugate $\overline{v}$.
By de Moivre's theorem, therefore, we find $f_a$ by replacing each $\cos(n t)$ with $\frac{1}{2}(v^n + \overline{v}^n)$ and each $\sin(n t)$ with $-\frac{\rmi}{2}(v^n - \overline{v}^n)$,
\begin{eqnarray}
\fl   f_{a}(u,v,\overline{v}) & = & 256 u^3 + 12 a^2 u (12 v^3 - 12 \overline{v}^3 + 22 v^2 - 22 \overline{v}^2 + 4 v  - 12 \overline{v}  - 5 ) \textstyle + \frac{1}{2} a^3 ( 384 + 27 v^5  \nonumber \\
\fl    & &  + 27 \overline{v}^5 + 44 v^4 + 172 \overline{v}^4 + 102 v^3 + 186 \overline{v}^3 + 460 v^2 + 52 \overline{v}^2 + 748 v - 4 \overline{v} ).      \label{eq:fa}
\end{eqnarray}
With the restriction $|v| = |\overline{v}| = 1$, $f_a$ coincides with $p_{t,a}$, and therefore in some sense its nodal set contains $5_2$.

The unit three-sphere $S^3 = \{ (u,v) \in \mathbb{C}^2; |u|^2 + |v|^2 = 1\}$ is contained in $\mathbb{C}^2$, i.e.~the domain of $f_a$, and we can therefore study the nodal set of $f_a$ restricted to $S^3$, for some appropriate $a$.
In general, we cannot expect this set, $f_a^{-1}(0) \cap S^3$, to be ambient isotopic to the closed braid $f_a^{-1}(0) \cap (\mathbb{C}\times S^1) \simeq p_{t,a}^{-1}(0)$, as these might be quite distant in $\mathbb{C}^2$.
However, the two sets approach as $a$ decreases to $0$, as follows.
The parameter $a$ determines the maximum absolute value any root of the polynomial $p_{t,a}(u)$ may have; from our previous discussion, it is $2.098 a$.
Thus, as $a$ decreases, the relevant range of $u$ is a smaller neighbourhood of $u$ in the complex $u$-plane; in this limit, the relevant $v$ in the 3-sphere (satisfying $|u|^2 + |v|^2 = 1$) are then of approximately unit modulus, coinciding with the braid space.
It ought to be sufficient, therefore, to find sufficiently small $a$.
As $a$ decreases, the nodal set of $f_a$ approaches the locus $u = 0$.

\begin{figure}[htbp]
\begin{center}
\includegraphics[width=14cm]{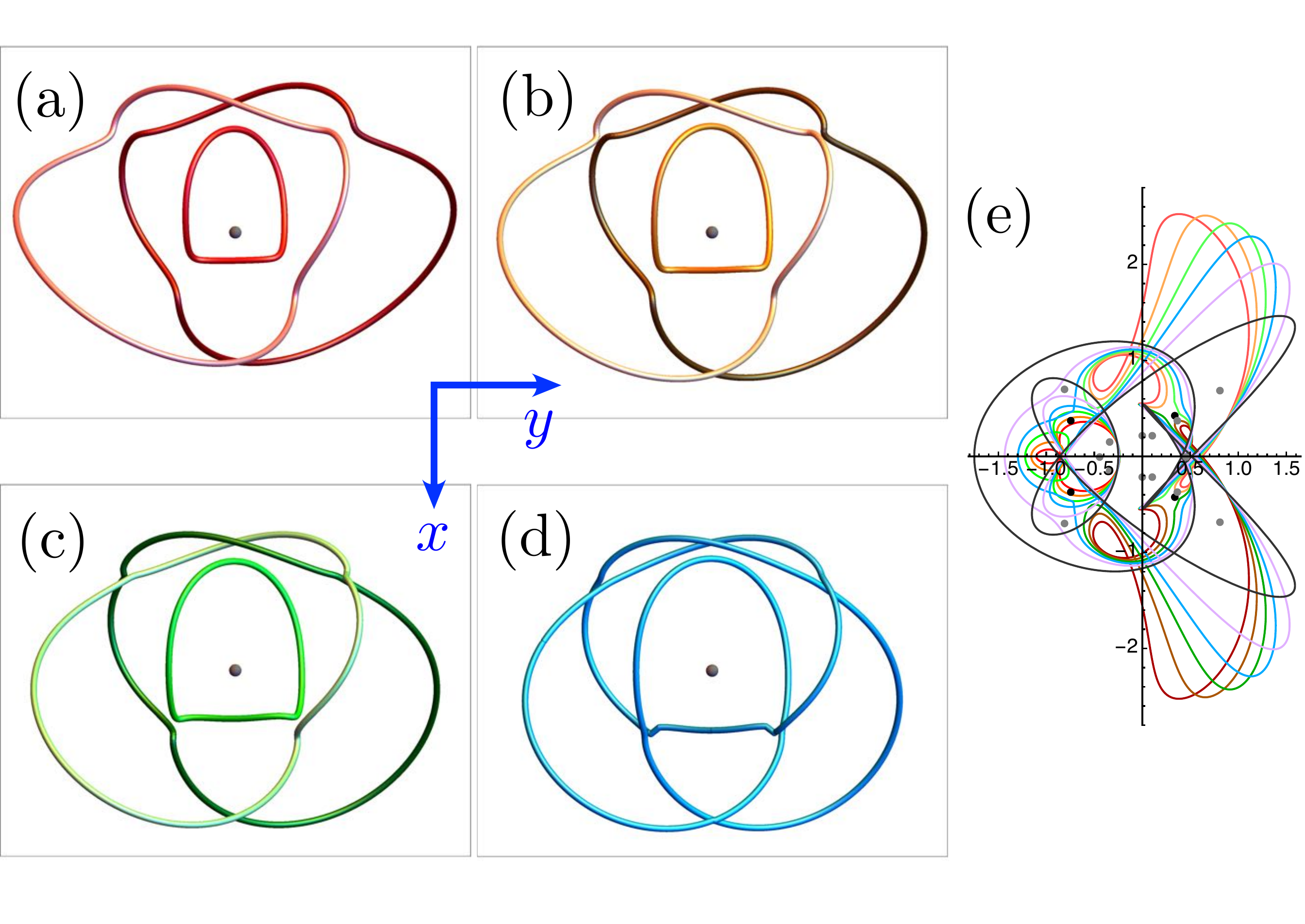}
\caption{Decreasing $a$ to the correct $5_2$ topology. 
   (a)-(d) represent the nodal set of $f_a(u,v,\overline{v})$ in stereographic coordinates (\ref{eq:stereo}), viewed down the $z$-axis from $(0,0,\infty)$: (a) $a = 0.7$; (b) $a = 0.6$; (c) $a = 0.5$; (d) $a = 0.4$.
   The sequence of reconnections leading to the $5_2$ knot in (d) is explained in the text.
   The origin is represented by a grey sphere.
   (e) Sequence of rescaled trajectories $a^{-1}(R(\phi),z(\phi))$ and coloured as (a)-(d), together with $a = 1/4$ (violet) and the limiting trajectory $T$ (black).
   Reconnection events, represented by black dots, occur at $a = 0.537\ldots$ and $a = 0.430\ldots$.
   $(R,z)$ coordinates of all the critical points of $|p_{t,a}|^2$ are shown as grey dots.
}
\label{fig:reconnect}
\end{center}
\end{figure}

In visualising the nodal sets for varying $a$, we use the stereographic representation (\ref{eq:stereo}) of $u$ and $v$, so that in cylindrical polar coordinates, $\phi$ corresponds to the braid parameter $t$.
Choosing $a = 1$ does not give $5_2$; its nodal set $f_{a=1}^{-1}(0) \cap S^3$ consists of three unknotted rings, one close to the $z$-axis near $z = 0$, and two larger rings around the $z$ axis, symmetrically arranged around the $z = 0$ plane (one above, one below).
This is also the case on decreasing $a$ by a small amount; Figure \ref{fig:reconnect} (a) shows the nodal lines for $a = 0.7$, viewed down the $z$-axis from $(0,0,+\infty)$.
From the way $f_a$ was constructed from the braid, we expect to be able to read the braid word from projection, with $\sigma_1$ occurring at larger radius than $\sigma_2$, and at a positive crossing as $\phi$ increases, the upper strand has increasing $R$. 
At $a = 0.7$, the word of this ring configuration is clearly $\sigma_1^{-1} \sigma_1$, consistent with an unlink/unknot.
Figure \ref{fig:reconnect} (b) shows, the nodal set for $a = 0.6$; although the topology has not changed, the two outer rings are approaching at two points at either side of the $\sigma_1$ crossing, locally resembling hyperbolae.
By $a = 0.5$, as in Figure \ref{fig:reconnect} (c), the outer rings have `reconnected' at the points where they had previously approached so the topology of the configuration has changed, with the braid word now $\sigma_1^{-1} \sigma_1^3$.
This is not yet $5_2$, but all of the $\sigma_1$ crossings are correct for $w$.
Towards the bottom of Figure \ref{fig:reconnect} (c), there are two new hyperbolic approaches on each side of the $\sigma_1^{-1}$ crossing.
Figure \ref{fig:reconnect} (d) shows the nodal set at $a = 0.4$, where there have been two further reconnections at the previous approaches, and the word on projection is indeed now $w$, meaning that this is a nodal $5_2$ knot.

This process of nodal line topology changing with $a$ can be visualized in terms of the trajectories of the nodal lines for different $a$ projected in the $(R,z)$ plane with $\phi$ as a parameter.
Figure \ref{fig:reconnect} (e) shows these curves, for various choices of $a$, appropriately rescaled to approach the braid trajectory $(X(t),Y(t))$ as $a \to 0$.
The nodal curve topology changes at two values of $a$, at each of which two reconnection events occur.
The nodal lines of a typical function like $f_a$, under evolution of a parameter $a$, reconnect at points where the nodal lines coincide with a critical point of absolute value $|f_a|^2$ (i.e.~isolated points where $\nabla |f_a|^2 = 0$) \cite{bd}.
Numerically solving the equation for $a$, $\nabla |f_a|^2 = 0$ coincides with $f_a = 0$ at $a = 0.537\ldots$, at which reconnections occur at $R = 1.184\ldots$, $\phi = \pi \pm 0.841\ldots$ and $z = \pm0.228\ldots$.
The other reconnection occurs at $a = 0.430\ldots$, with reconnection coordinates $R = 0.680\ldots, \phi = \pm 0.827\ldots, z = \pm 0.158\ldots.$
The $(R,z)$ coordinates of these 14 reconnection points are shown as black dots in Figure \ref{fig:reconnect} (e); as the trajectories of nearby $a$ approach these points, they undergo hyperbolic reconnection.

Although we do not prove that there are no further reconnections for $a < 0.4$, the following argument gives some evidence that this does not happen.
As $a \to 0$, in a solid torus neighbourhood of the unit circle in the $z = 0$ plane, $f_a$ approaches $p_{t,a}(u)$.
Since $p_{a,t}(u)$ is holomorphic in $u$, critical points of argument lie on the two families of complex saddles $c_{\pm}(t)$ where $\rmd p_{a,t}/\rmd u = 0$.
All of the critical points of $|p_{t,a}(u)|^2$ can be located by finding the extrema of $p_{t,a}(c_{\pm}(t))$ with respect to $t$, and these are shown as grey dots in Figure \ref{fig:reconnect} (e).
For any finite $a$, $f_a$ will have critical points corresponding to these (albeit at slightly different coordinates), as can be seen for the reconnection events which do occur.
The trajectory at $a = 0.4$, and the limiting $(X,Y)$ trajectory, do not approach any other critical points, so we do not expect further reconnections (nor do we expect further critical points to occur, although we have not shown this).

\begin{figure}
\begin{center}
\includegraphics[width=12.6cm]{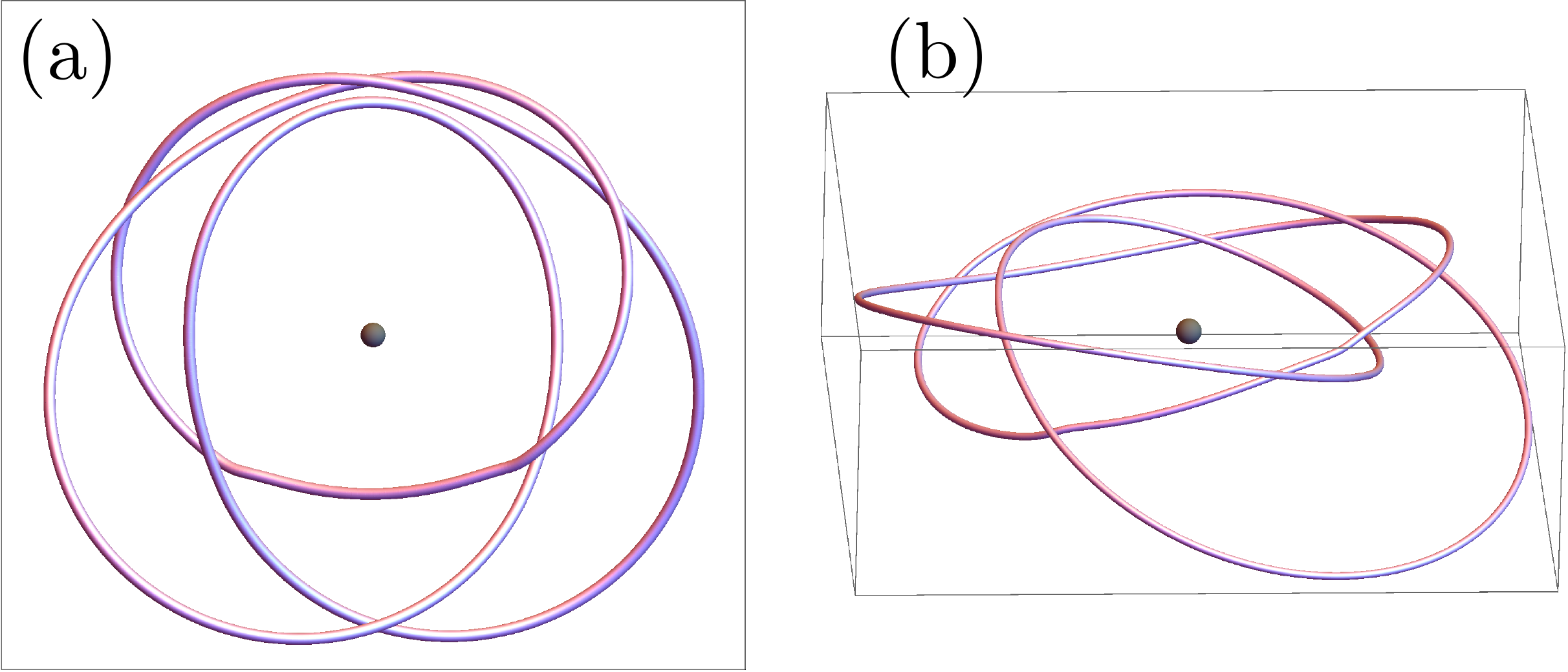}
\caption{The nodal set, a $5_2$ knot, of the function (\ref{eq:f5_2}) with $a = \frac{1}{4}$, in stereographic coordinates (\ref{eq:stereo}).
(a) vertical view; (b) side view.
As in Figure \ref{fig:reconnect}, the grey sphere is the origin.
}
\label{fig:5_2}
\end{center}
\end{figure}

We therefore have constructed a family of functions $f_{5_2}$ whose nodal set is the three-twist knot $5_2$: $f_a$ of (\ref{eq:fa}), restricted to the 3-sphere with any $a < 0.430\ldots$ will suffice.
With the choice of explicit coordinates (\ref{eq:stereo}) and $a = \frac{1}{4}$, this is
\begin{eqnarray}
\fl   f_{5_2} & = & 2042 x^{10} + 10210 x^8 y^2  + 10210 x^8 z^2 + 20420 x^6 y^4 + 40840 x^6 y^2 z^2 + 20420 x^6 z^4 \nonumber \\
\fl   & & + 20420 x^4 y^6 + 20420 x^4 z^6 + 61260 x^4 y^4 z^2 + 61260 x^4 y^2 z^4 + 10210 x^2 y^8 + 40840 x^2 y^6 z^2 \nonumber \\
\fl   & & + 61260 x^2 y^4 z^4 + 40840 x^2 y^2 z^6 + 10210 x^2 z^8 + 2042 y^{10} + 10210 y^8 z^2 + 20420 y^6 z^4 \nonumber \\
\fl   & & + 20420 y^4 z^6 + 10210 y^2 z^8 + 2042 z^{10} - 3 x^9 - 12 x^7 y^2 - 12 x^7 z^2 - 18 x^5 y^4 - 36 x^5 y^2 z^2 \nonumber \\
\fl   & & - 18 x^5 z^4 - 12 x^3 y^6 - 36 x^3 y^4 z^2 - 36 x^3 y^2 z^4 - 12 x^3 z^6 - 8328 x^2 y^6 - 3 x y^8 - 12 x y^6 z^2 \nonumber \\
\fl   & & - 12 x y^2 z^6 - 18 x y^4 z^4 - 3 x z^8 - 1890 x^8 - 7816 x^6 y^2 - 384 x^6 y z - 32264 x^6 z^2 \nonumber \\ 
\fl   & & - 12108 x^4 y^4 - 1152 x^4 y^3 z - 97560 x^4 y^2 z^2 - 1152 x^4 y z^3 - 85452 x^4 z^4 - 1152 x^2 y^5 z \nonumber \\
\fl   & & - 98328 x^2 y^4 z^2 - 2304 x^2 y^3 z^3 - 171672 x^2 y^2 z^4 - 1152 x^2 y z^5 - 81672 x^2 z^6 - 2146 y^8 \nonumber \\
\fl   & & - 384 y^7 z - 33032 y^6 z^2 - 1152 y^5 z^3 - 86220 y^4 z^4 - 1152 y^3 z^5 - 81928 y^2 z^6 - 384 y z^7 \nonumber \\
\fl   & & - 26594 z^8 + 324 x^7 + 396 x^5 y^2 - 4224 x^5 y z + 828 x^5 z^2 - 180 x^3 y^4 - 8448 x^3 y^3 z \nonumber \\
\fl   & & + 504 x^3 y^2 z^2 - 8448 x^3 y z^3 + 684 x^3 z^4 - 252 x y^6 - 4224 x y^5 z - 324 x y^4 z^2 - 8448 x y^3 z^3 \nonumber \\
\fl   & & + 108 x y^2 z^4 - 4224 x y z^5 + 180 x z^6 - 3316 x^6 - 12444 x^4 y^2 - 8064 x^4 y z - 35340 x^4 z^2 \nonumber \\
\fl   & & - 13212 x^2 y^4 - 6912 x^2 y^3 z - 73944 x^2 y^2 z^2 - 9216 x^2 y z^3 - 60516 x^2 z^4 - 4084 y^6 \nonumber \\
\fl   & & + 1152 y^5 z - 36876 y^4 z^2 - 61284 y^2 z^4 - 1152 y z^5 - 28492 z^6 + 954 x^5 - 540 x^3 y^2 \nonumber \\
\fl   & & - 8448 x^3 y z + 1404 x^3 z^2 + 234 x y^4 + 252 x y^2 z^2 - 8448 x y z^3 + 558 x z^4 + 4996 x^4 \nonumber \\
\fl   & & + 7496 x^2 y^2 - 8064 x^2 y z + 33752 x^2 z^2 - 8448 x y^3 z + 4228 y^4 + 1152 y^3 z + 32984 y^2 z^2 \nonumber \\
\fl   & & - 1152 y z^3 + 28972 z^4 + 708 x^3 + 132 x y^2 - 4224 x y z + 564 x z^2 + 2386 x^2 + 2130 y^2 \nonumber \\
\fl   & & - 384 y z + 26834 z^2 + 189 x - 1994 + 2\rmi \left( 1+ x^2 + y^2 + z^2\right)^3\left( 143 x^6 y + 6114 x^6 z \right.\nonumber \\
\fl   & & + 429 x^4 y^3 + 18342 x^4 y^2 z + 429 x^4 y z^2 + 18342 x^4 z^3 + 429 x^2 y^5 + 18342 x^2 y^4 z \nonumber \\
\fl   & & + 858 x^2 y^3 z^2 + 36684 x^2 y^2 z^3 + 429 x^2 y z^4 + 18342 x^2 z^5 + 143 y^7 + 6114 y^6 z + 429 y^5 z^2 \nonumber \\
\fl   & & + 18342 y^4 z^3 + 429 y^3 z^4 + 18342 y^2 z^5 + 143 y z^6 + 6114 z^7 + 1158 x^5 y - 96 x^5 z \nonumber \\
\fl   & & + 2316 x^3 y^3 - 192 x^3 y^2 z + 2316 x^3 y z^2 - 192 x^3 z^3 + 1158 x y^5 - 96 x y^4 z +2316 x y^3 z^2 \nonumber \\
\fl   & & -192 x y^2 z^3 + 1158 x y z^4 -96 x z^5 + 1902 x^4 y - 6234 x^4 z + 1584 x^2 y^3 - 12468 x^2 y^2 z \nonumber \\
\fl   & & + 2139 x^2 y z^2 - 20660 x^2 z^3 - 318 y^5 - 6234 y^4 z - 81 y^3 z^2 - 20660 y^2 z^3 + 237 y z^4 \nonumber \\
\fl   & & - 14426 z^5 - 52 x^3 y - 192 x^3 z + 460 x y^3 - 192 x y^2 z + 204 x y z^2 - 192 x z^3 - 1746 x^2 y \nonumber \\
\fl   & & - 6234 x^2 z + 642 y^3 - 6234 y^2 z + 45 y z^2 - 14426 z^3 - 954 x y - 96 x z - 49 y + 6114 z \nonumber \\
\fl   & & \left. + 6114 z \right) 
      \label{eq:f5_2}
\end{eqnarray}
where (\ref{eq:fa}) has been multiplied by $8$ to ensure the coefficients are Gaussian integers.
The nodal knot of (\ref{eq:f5_2}) is shown in Figure \ref{fig:5_2}, where as with Figure \ref{fig:reconnect} (d) the vertical projection enables the braid word $w$ to be clearly read off as $\phi$ increases from $0$ to $2\pi$.
This explicit function $f_{5_2}(x,y,z)$ could be embedded directly in physical systems as summarized in Section \ref{sec:int}.
Although it is not necessarily obvious from the Figures, as $a$ decreases, the nodal lines of $f_a$ get closer to the unit circle in the $z = 0$ plane; suitably large values of $a$ are required to ensure the nodal set is the appropriate knot by inspection, as in the examples here.

\section{Critical points of argument and Morse-Novikov number}\label{Sec:mn}

\begin{figure}
\begin{center}
\includegraphics[width=12.87cm]{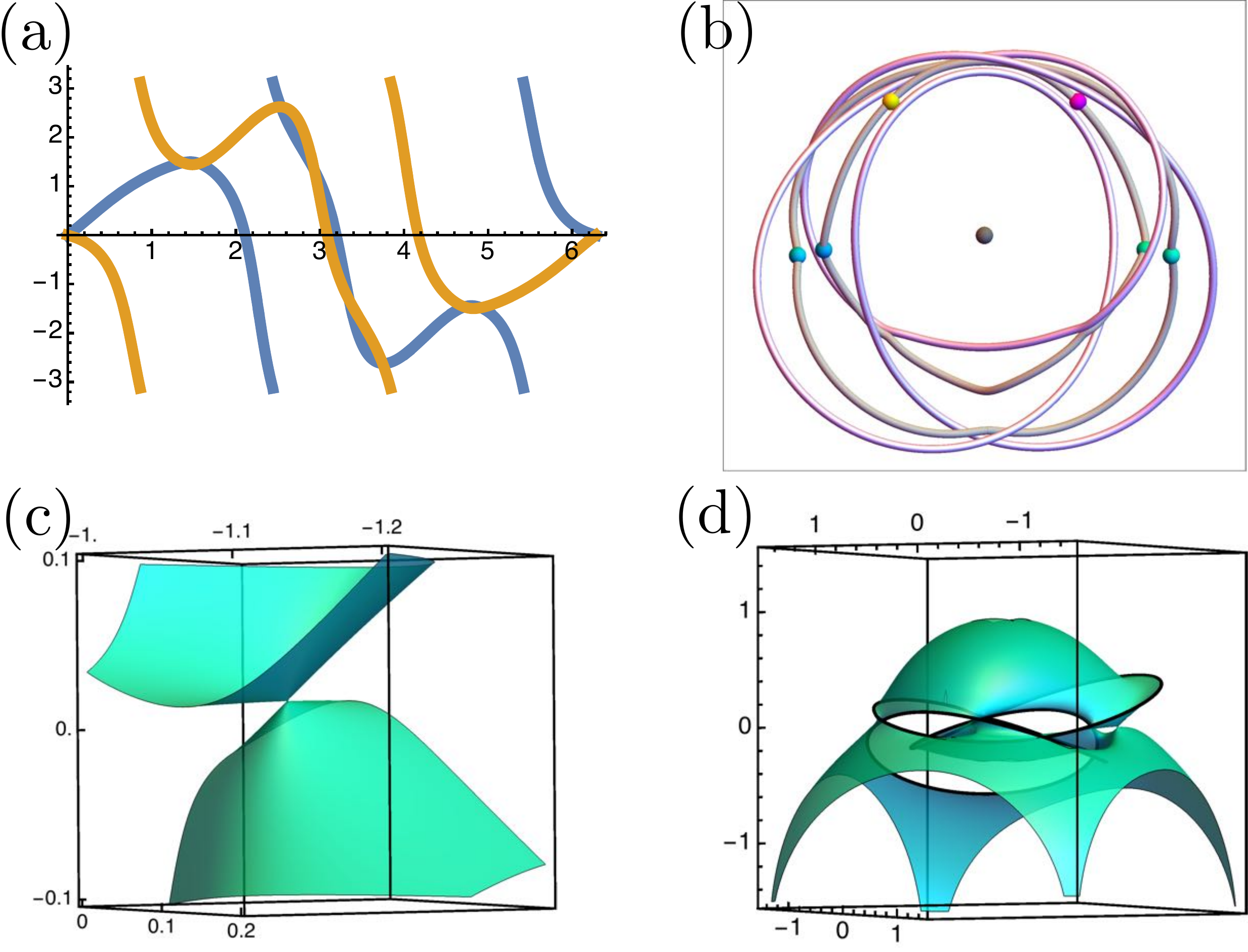}
\caption{
   The nonfibred knot $5_2$ has critical points of complex argument.
   (a) Plots of the argument on the critical lines of the polynomial, $\arg(p_{t,a}(c_{\pm}(t)))$ as functions of $t$.
   There are six maxima and minima.
   (b) $5_2$ nodal knot for $a = 1/4$.
   The positions of the critical points are given by the coloured spheres (with hue corresponding to values of argument at the critical points), which are numerically very close to the critical values of the polynomial in (a).
   The grey curve here is the locus of the critical values where $\nabla_{R,z}\arg f_{5_2} = 0$ for each $\phi$, i.e.~transverse argument saddle points; this curve is close to the critical loci of the polynomial $c_{\pm}(t)$.
   (c) Level set of constant argument $1.499$, corresponding to the first argument critical point.
   (d) The Seifert level set (excluding the critical point) of this value of argument.  
}
\label{fig:nonfibred}
\end{center}
\end{figure}

With this the explicit $f_{5_2}$ function we now explore how $\arg{f_{5_2}}:S^{3}\backslash 5_{2}\longrightarrow S^{1}$ fails to be a fibration. 
As discussed above, since the three-twist knot $5_2$ is not fibred, any $f_{5_2}$ must have critical points of phase (i.e.~argument), where $\nabla \arg f_{5_{2}}=0$. 
Following the discussions in \cite{lemniscate} and \cite{proof} we know that for sufficiently small $a$, the critical points of $\arg f_{a}$ in $S^3$ are well approximated by those of $\arg p_{t,a}(u)$ as a function of $u$ and $t$, if $\mathbb{C}\times S^{1}$ is embedded into $\mathbb{C}^{2}$ in the usual way by $(u,t)\mapsto(u,\rme^{\rmi\,t})$. 
Just as for critical points of absolute value, the critical points of the argument of $p_{t,a}(u)$ lie on the critical curves $c_{\pm}$ where $\rmd p_{t,a}(u)/\rmd u = 0$.
The arguments $\arg p_{t,a}(c_{\pm}(t))$ as functions of $t$ are shown in Figure \ref{fig:nonfibred} (a): evidently there are six maxima and minima, numerically located at $t = 1.46\ldots, 1.48\ldots, 2.52\ldots, 3.76\ldots, 4.80\ldots, 4.82\ldots$.
The critical points of $f_{5_2}$ in (\ref{eq:f5_2}) when $a = \frac{1}{4}$ are correspondingly close, they are plotted around the nodal knot in Figure \ref{fig:nonfibred} (b).
The level set of constant argument at one of these phase critical points (the one with lowest value of $\phi$), with argument $1.499\ldots$, is shown in Figure \ref{fig:nonfibred} (c).
At the critical point this is a self-intersecting surface locally resembling a diabolo (double cone), and for nearby values of phase, the levels sets resemble hyperboloids of one or two sheets (not shown), and hence have different genus.
This constant phase level set fails to be a surface at the critical point; apart from this point it resembles a Seifert surface for the knot (Figure \ref{fig:nonfibred} (d)).

It can be shown that $\arg(f_{a})$ is a smooth, circle-valued Morse function on the knot complement $S^{3}\backslash 5_{2}$. 
It is also \emph{regular}, that is in a solid torus neigbourhood of the knot, each level set of phase is a smooth ribbon.
The minimal number $m$ of critical points of any smooth, circle-valued, regular Morse function on $S^3\backslash K$ is zero if and only if the knot $K$ is fibred; this is the Morse-Novikov number of the knot $K$.
Since $5_{2}$ is not fibred, $m$ cannot be zero and in fact the Morse-Novikov number of $5_{2}$ is known to be $m = 2$ \cite{goda}. 
The argument of the function $f_{a}$ that we constructed has six critical points and hence does not actually have the minimal number of critical points.
However, in principle the two nearby maximum-minimum pairs in Figure \ref{fig:nonfibred} (a) could be smoothed away, leaving only two; indeed in \cite{proof} we show that there exists a more complicated, finite Fourier parametrization of the same braid word that only has two phase critical points.

\section{How $f_{5_{2}}$ fails to have $5_2$ as the link of a weakly isolated singularity}\label{sec:weakly}

The function $f_{5_{2}}$ was constructed to have a nodal set on the unit three-sphere which is the three-twist knot $5_{2}$.
As such, it was defined as the restriction of the function $f_a$ to $S^3$ for sufficiently small, fixed $a$; the domain of $f_a$ is a four-dimensional space, which may be considered as $\mathbb{C}^2$ (with coordinates $(u,v)$), $\mathbb{C}\times\mathbb{R}^2$ (as $f_a$ is semiholomorphic), or $\mathbb{R}^4$; we will use $\mathbb{R}^4$ in the following discussion.

The topology of the nodal sets of $f_a$ on three-spheres around $(0,0,0,0)\in\mathbb{R}^{4}$ of radius other than unity might not be $5_{2}$.
This clearly distinguishes our construction from the more general constructions of Brauner \cite{brauner} for torus knots, and Perron \cite{perron} and Rudolph \cite{rudolph} for the figure-eight knot, all of whose functions $f$ are proved to have the desired knot as their nodal set on all three-spheres around $(0,0,0,0)$ whose radius is sufficiently small. 
Furthermore, they all have an isolated singular point at the origin of $\mathbb{R}^4$, meaning that $f(0,0,0,0) = 0$, $\nabla f|_{(0,0,0,0)}$ does not have full rank at the origin and $\nabla f(x)$ has full rank for all other $x\in \mathbb{R}^4$ in a neighbourhood of the origin.
A theorem by Milnor states every knot which arises in this way must be fibred \cite{milnor}. 
Since $5_2$ is not fibred, such a function for $5_2$ does not exist.

A weaker condition is for the the singular point at the origin to be only \emph{weakly isolated}, that is, we demand $f(0,0,0,0) = 0$, $\nabla f|_{(0,0,0,0)}$ does not have full rank and $\nabla f(x)$ has full rank only for all $x$ in a neighbourhood of the origin that satisfy $f(x) = 0$.
Akbulut and King \cite{ak} showed that any knot can arise as the link around a weakly isolated singular point of a polynomial $\mathbb{R}^{4}\rightarrow\mathbb{R}^{2}$. 
However, the polynomial $f_a$ given by (\ref{eq:fa}) is clearly not of this form, since $(0,0,0,0)$ is not even part of its nodal set.

The origin is, however, part of the nodal set of modifications of $f_a$ of the form $F = (v\overline{v})^{ks}f_{a}\left(u/(v\overline{v})^{k},v\right)$ for all natural number $k$; the origin is also a singular point and can easily be shown to be weakly isolated.
It can also be shown that for small $a$ the nodal set of $F$ on the unit-three sphere is still the knot $5_{2}$. 
Thus, if for all $\rho\in(0,1]$ the value $0$ is a regular value of $F$ restricted to $S^{3}_{\rho}$ (the three-sphere of radius $\rho$), then the nodal set of $F$ on $S^{3}_{\rho}$ is indeed the knot $5_{2}$ and $F$ is hence of the form discussed in \cite{ak}.
For small $a$ and large $k$, the intersection $F^{-1}(0) \cap S^{3}_{\rho}$ at $(u,v)$ is transverse for all $\rho\in(0,1]$ if $u$ is a simple root of the complex polynomial $F(\bullet,v)$. 
Thus $0$ is a regular value of $F|_{S^{3}_{\rho}}$ for all $\rho\in(0,1]$ if for all $r \in (0,1]$ and all $t \in [0,2\pi)$ the polynomial $F(\bullet,r \rme^{\rmi t})$ has only simple roots. 
This is equivalent to the statement that $f_{5_{2}}(\bullet,r \rme^{\rmi t})$ has only simple roots for all $r\in(0,1]$ and all $t\in[0,2\pi)$.

However, for our polynomial $f_a$ in (\ref{eq:fa}), this is not the case. 
The values of $r$ and $t$ for which $f_{5_2}$ has non-simple zeros are identical to the zeros of the resultant of the polynomials $f_{a}$ and its derivative $\rmd f_a/\rmd u$. 
We find four such points for $(r,t)$. 
The corresponding polynomials $F(r\rme^{\rmi t})$ have non-simple roots and at these points zero is not a regular value of $F|S_{\rho}^{3}$.
In fact, the discussion in Section \ref{sec:embed} of reconnection events on the nodal set of $f_{a}$ on $S^{3}$ as $a$ varies also tells us that these non-simple roots exist. 
The points $(u,r \rme^{\rmi t})\in\mathbb{C}^{2}$ where the reconnections occur are double roots of $f_{a}(\bullet,re^{it})$. 
Since $a$ is only a scaling parameter, this means that $f_{a}(\bullet,r \rme^{\rmi t})$ has a double root for any positive value of $a$.
We find that on the corresponding points the gradient $\nabla_{S_{\rho}^{3}} F$ does not have full rank. 
Hence the function $F$ has a weakly isolated singularity at $(0,0,0,0)$ and the knot $5_{2}$ as its zero level set on the unit three-sphere, but its nodal sets on three-spheres of small radius are not $5_{2}$.

\section{Polynomials whose nodal set are other knots of six and fewer crossings}\label{sec:examples}

Our focus in this paper has been finding the complex function $f_{5_2}$ for the three-twist knot $5_2$, the simplest knot not previously found as a nodal knot for an explicit complex map.
However, the procedure can be easily extended to any knot or link via chosen braid representation; we now summarize examples of knots and links with 6 and fewer crossings not described by previous discussions.

As described in the introduction, all of the knots of five or fewer crossings except $5_2$ have been described previously, as they are either torus knots (and hence functions constructed by Brauner's method \cite{brauner}) or $4_1$ which is a lemniscate knot \cite{isolated,lemniscate}.
The 6-crossing knots include the prime knots $6_1$ (the stevedore's knot), $6_2$ and $6_3$, and the two composite $3_1 \# 3_1$ (the granny knot) and $3_1 \# 3_1^*$ (the square knot).
In fact, $6_3$ is also a lemniscate knot (the $(5,2,2)$ lemniscate knot) \cite{lemniscate}, and the granny knot can be constructed by a generalization of the lemniscate construction \cite{lemniscate}, so we focus here on $6_2$, $6_1$ and $3_1 \# 3_1^*$.

\begin{figure}
\begin{center}
\includegraphics[width=14.2cm]{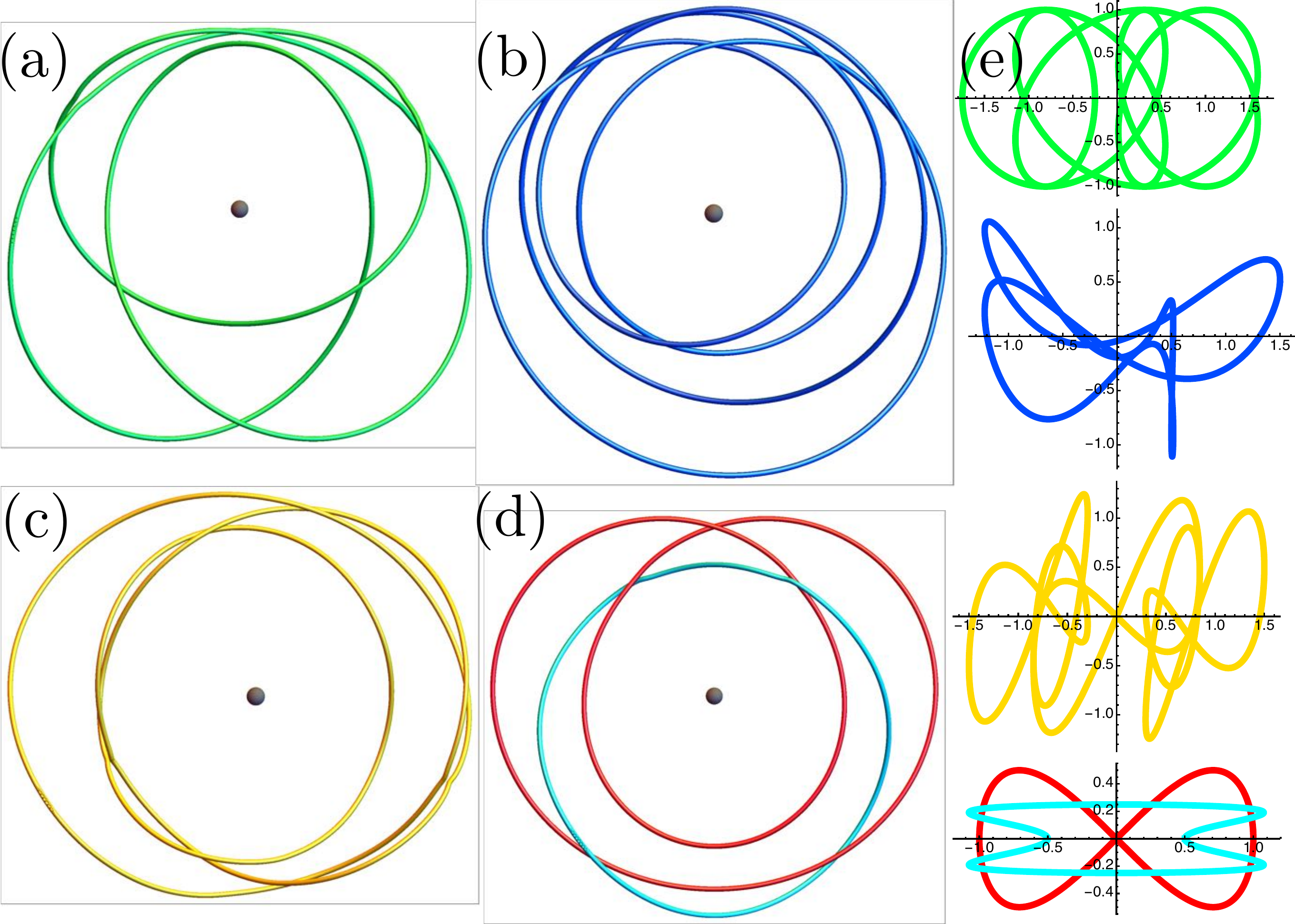}
\caption{Nodal sets of functions for six-crossing knots and Whitehead link.
(a) $6_2$ knot; (b) $6_1$ knot; (c) $3_1\#3_1^*$ composite knot; (d) Whitehead link $L5a1$.
(e) The respective braid trajectories used to generate the functions with these nodal sets, given by equations (\ref{eq:T6_2}), (\ref{eq:T6_1}), (\ref{eq:Tsquare}) and (\ref{eq:whitehead}) respectively.
For each of the functions, the parameter value $a = \frac{1}{4}$ was chosen. 
}
\label{fig:extras}
\end{center}
\end{figure}

The minimal braid word for $6_2$ is $\sigma_1^3 \sigma_2 \sigma_1 \sigma_2$, which has exactly the same form as $w$ for $5_2$, except for the sign of one occurrence of $\sigma_1$.
Therefore a suitable trajectory $T_{6_2}$ for this knot can use (\ref{eq:Xchoice}) for $X(t)$, and in fact, a suitable $Y(t)$ which guarantees the correct signs of the braid word is the single Fourier term $-\sin(\frac{5}{3}t)$, so
\begin{equation}
\textstyle   T_{6_2} : \left( - \cos(\frac{2}{3}t) - \frac{3}{4} \cos(\frac{5}{3}t), -\sin(\frac{5}{3}t)\right).
   \label{eq:T6_2}
\end{equation}
We then follow exactly the same procedure as above, and find similar reconnection events as $f_{5_2}$ at similar values of $a$.
The trajectory and nodal set of $f_{6_2}$ at $a = \frac{1}{4}$ are shown in Figure \ref{fig:extras} (a).

The minimal braid word for the stevedore knot is more complicated, having $4$ strands: $\sigma_3^{-1}\sigma_1\sigma_2\sigma_1^{-1}\sigma_3^{-1}\sigma_2\sigma_1$ (which is a simple Markov conjugation away from the standard minimal form \cite{gittings}).
The corresponding unsigned word is rather symmetric (placing the adjacent $\sigma_1^{-1}\sigma_3^{-1}$ at the same value of $t$), and a suitable function $Y(t)$ is $\cos(\frac{1}{2}t) + \frac{1}{2} \cos(\frac{3}{4} t)$.
We could not find any two-term Fourier series as an $Y(t)$ which gave the correct crossings, but a numerical search amongst three-term series was successful, albeit introducing rational angle phase shifts; our resulting trajectory for $6_1$ is
\begin{equation}
\fl \textstyle   T_{6_1} : \left( \cos(\frac{1}{2}t) + \frac{1}{2} \cos(\frac{3}{4} t), \cos(\frac{1}{4}[t+\pi]) - \frac{2}{5} \cos(\frac{1}{8}[4t+ \pi]) + \cos(\frac{1}{4}[5t+\pi])\right).
   \label{eq:T6_1}
\end{equation}
The phase shifts mean that the corresponding Fourier series (and polynomial) has irrational coefficients.
Despite this extra complication, the rest of the construction follows through in the same way as the others, and $6_1$ is the nodal set for similar values of $a$ to the knots previously considered.
The trajectory and nodal set of $f_{6_1}$ at $a = \frac{1}{4}$ are shown in Figure \ref{fig:extras} (b).

The square knot is one of the simplest composite knots, consisting of two trefoils with opposite handedness.
It has minimum word $\sigma_1^3 \sigma_2^{-3}$; the granny knot consists of two trefoils of the same handedness, with word $\sigma_1^3 \sigma_2^{3}$.
Functions for which these are the nodal set can be constructed from a braid function with the same $Y(t)$ function, and a choice is $\sin(\frac{1}{3}t) + \frac{1}{2} \sin(\frac{5}{3}t)$.
A suitable crossing function $Y(t)$ for the granny knot is $-\cos(\frac{5}{3}t)$.
Finding the analogous function for the square knot, whose crossing signs are different, is a bit complicated, but can be achieved with a two-term Fourier series, with the resulting trajectory
\begin{equation}
\textstyle   T_{3_1\#3_1^*} : \left(\sin(\frac{1}{3}t) + \frac{1}{2} \sin(\frac{5}{3}t),\frac{1}{2} \sin(\frac{5}{3}t) + \frac{3}{4} \sin(\frac{8}{3}t)\right).
   \label{eq:Tsquare}
\end{equation}
The trajectory and nodal set of $f_{6_1}$ at $a = \frac{1}{4}$ are shown in Figure \ref{fig:extras} (c).

We conclude this section by considering an example of a link, namely the Whitehead link ($L5a1$ in the Thistlethwaite table \cite{knotatlas,knotinfo}, $5_1^2$ in the Rolfsen table \cite{rolfsen}).
The procedure for links follows the same argument as for knots; each component of the link is comprised of a subset of the braid's strands.
The minimal braid for the Whitehead link has five crossings and three strands, with minimal word is $\sigma_1 \sigma_2 \sigma_1 \sigma_2 \sigma_1 \sigma_2$.
Here, strands $2$ and $3$ are the same component, with a single occurrence of their generator: the trajectory of this component is a simple loop.
The other loop winds around this one, in a lemniscate-like trajectory. 
A choice of the trajectory pair making the Whitehead link is therefore the following:
\begin{equation}
   T_{L5a1} : \left\{ \begin{array}{l} \left( \cos(\frac{1}{2}t) - \frac{1}{2} \cos(\frac{3}{2}t), \frac{1}{4} \sin(\frac{1}{2}t) \right), \\ 
   \left( \cos(t), \frac{1}{2}\sin(2t)\right) \end{array} \right.
   \label{eq:whitehead}
\end{equation}
The trajectory and nodal set of $f_{L5a1}$ at $a = \frac{1}{4}$ are shown in Figure \ref{fig:extras} (d).
A simple way of creating functions with nodal links is by multiplying functions whose nodal set is each component loop of the link, so the braid construction for knots is less necessary for constructing functions with nodal links.
Not all link functions will factorize in this way, and in fact the example $f_{L5a1}$ here does not.
There are six 6-crossing prime links, of which two are lemniscate ($L6a3$ is the torus link $T_{6,3}$ and $L6a4$ is the borromean rings); braid functions for the others are straightforward (albeit tedious) to derive, and we omit them here.

\section{Concluding remarks}\label{conc}

We have described a construction of a function from a the unit three-sphere or three-dimensional cartesian space to the complex numbers whose nodal set is the three-twist knot $5_2$, with an approach based on a braid representative that can be extended (with sufficient patience) to any knot or link.
The same $X(t)$ function was used with a different $Y(t)$ to make the $6_2$ knot, but other generalizations are possible; a simple further case is to take two copies of the braid word before closing, equivalent to replacing every occurrence of $v$ and $\overline{v}$ in (\ref{eq:fa}) with $v^2$ and $\overline{v}^2$.
The knot corresponding to the closure of $w^2$ is the 12-crossing knot $12n_{0750}$ \cite{knotinfo}.

The procedure described here to find the correct braid trajectory was focused on getting the simplest possible Fourier series (based on number of terms, rational simplicity of coefficients and low orders of spatial frequencies).
The resulting polynomial (\ref{eq:fa}) was comparatively small, and had integer coefficients; these properties might be desirable, but are not necessary.
In particular, the nodal knot will be structurally stable to small perturbations of the coefficients, and indeed such a perturbation might lead to a conformation of the knot with better-separated strands \cite{isolated}.
It is not immediately obvious whether more careful choices of the Fourier series for the trajectory may give a simpler polynomial function as an alternative to $f_a$; as a polynomial, the order in the holomorphic parameter is equal to the number of strands $s$ of the braid representation, and the order of $v$, $\overline{v}$ is the highest order of $t/s$ is the trajectory Fourier series.
We expect the maximum value of the parameter $a$ at which the knot appears to depend on the order of the Fourier series, as well as the maximum distance of the trajectory from the origin; for instance, for the alternative $Y(t)$ function mentioned at the end of Section \ref{sec:braid}, the knot only appeared for $a \lesssim 0.13$, somewhat smaller than the function considered here (and therefore with the knot correspondingly closer to the unit circle).

It may be more practical for more complicated to make the process of determining the braid word more algorithmic, and in \cite{proof} we describe such an algorithm based on trigonometric interpolation which can be applied a braid word for any knot or link; this algorithm, applied to $5_2$, gives Fourier series for $X(t)$ and $Y(t)$ with real valued non-integer coefficients of six and eleven terms respectively, going up to order $8t/3$.

Being a complex scalar functions, it is natural to ask whether $f_{5_2}(x,y,z)$ might be realised in a scalar quantum wavefunction.
Indeed, despite not being normalizable, one could take $f_{5_2}(x,y,z)$ and evolve it in time according to the free Schr\"odinger equation directly (the time-dependent wavefunction also being a polynomial).
As with other similar constructions, however, the nodal set evolves, undergoing several reconnections after which the knot dissolves.
It is less easy to identify a system (e.g.~specified by a potential) in which a wavefunction involving $f_{5_2}$ is a time-independent eigenfunction.
Although it is possible to constructed certain torus links in nodal sets of eigenfunctions of the hydrogen atom \cite{berry}, most functions of Brauner form are incompatible with the symmetries of the potential \cite{taylor}.
It may appear more hopeful to consider $5_2$ as optical vortices in a propagating laser beam, extending \cite{isolated}, by taking the polynomial $f_{5_2}|_{z=0}$ in (\ref{eq:f5_2}) as an initial condition (hologram plane) for paraxial evolution \cite{polynomial}.
Preliminary investigation, however, fails to give an optical vortex $5_2$ knot for any value of $a$, and ongoing investigations are attempting to understanding this more deeply.

It should be possible to embed $f_{5_2}$ (or a perturbation) into the initial conditions for many other physical systems as discussed in Section \ref{sec:int}, such as (super)fluids, reaction-diffusion systems, Skyrme-Faddeev topological solitons and liquid crystals.
The resulting dynamics of the resulting knotted fields would be more complicated than simple linear Schr\"odinger evolution; it would be interesting to study these systems in more detail, to understand whether the evolution of different knots such as $5_2$ have special physics, possibly related to knot-theoretic properties such as being fibred.

\ack We are grateful to Gareth Alexander, Michael Berry, David Foster, Thomas Machon, Antoine Remond-Tiedrez, Jonathan Robbins, Danica Sugic, Alexander Taylor and Stu Whittington for discussions.
This work is funded by the Leverhulme Trust Research Programme Grant RP2013-K-009, SPOCK: Scientific Properties Of Complex Knots.

\end{document}